\documentclass[10pt]{article}
\usepackage[left=3cm,top=3cm,right=3cm,bottom=3cm]{geometry}
\usepackage{amsfonts,amssymb,amsmath,latexsym,indentfirst}
\usepackage{fontenc}
\usepackage{latexsym,wasysym,mathrsfs,enumerate}
\numberwithin{equation}{section}

\newtheorem{remark}{Remark}[section]
\newtheorem{definition}{Definition}[section]
\newtheorem{lemma}{Lemma}[section]
\newtheorem{theorem}{Theorem}[section]

\newcommand{\R}{\ensuremath{{\mathbb{R}} }}

\newcommand{\rr}{{\mathbb R}}

\newcommand{\ff}{\varphi}

\newcommand{\what}{\widehat}
\newcommand{\dd}{{\rm d}}
\newcommand{\ee}{{\rm e}}
\newcommand{\ii}{{\rm i}}

\newcommand{\peq}{\hspace*{0.10in}}
\newcommand{\ppeq}{\hspace*{0.05in}}
\newcommand{\para}{\hspace*{0.25in}}

\newcommand{\fim}{\hfill$\square$\\ \\}

\newcommand{\proof}{\noindent\textbf{Proof.}\quad}

\author{{\bf Amin Esfahani}\\  {\small School of Mathematics and Computer Science} \\
{\small Damghan University}\\
{\small  Damghan, Postal Code 36716-41167, Iran} \\\vspace{2mm}
{\small  E-mail: amin@impa.br, esfahani@du.ac.ir}\\
{\bf Luiz Gustavo Farah\footnote{Supported by FAPEMIG-Brazil and
CNPq-Brazil.}}\\
{\small ICEx, Universidade Federal de Minas Gerais}\\
{\small  Av. Ant\^{o}nio Carlos, 6627, Caixa Postal 702, 30123-970} \\
{\small  Belo Horizonte-MG, Brazil}\\
{\small  E-mail: lgfarah@gmail.com}
 }

\title{Local well-posedness for the Sixth-Order Boussinesq Equation}
\date{}

\begin{document}
\maketitle

\begin{abstract}
This work studies the local well-posedness of the initial-value problem for the nonlinear sixth-order
Boussinesq equation $u_{tt}=u_{xx}+\beta u_{xxxx}+u_{xxxxxx}+(u^2)_{xx}$, where $\beta=\pm1$. We prove local well-posedness
with initial data in non-homogeneous Sobolev spaces $H^s(\R)$ for negative indices of $s \in \R$.
\end{abstract}

\section{Introduction}
The study of wave propagation on the surface of water has been a subject of considerable theoretical and practical importance during the past decades. In 1872, Joseph Boussinesq \cite{boussinesq} derived a model equation for propagation of water waves from Euler's equations of motion for two-dimensional potential flow beneath a free surface by introducing appropriate approximations for small amplitude long waves.  Later the Boussinesq equation
\begin{equation}\label{fourthB}
u_{tt}=u_{xx}+\beta u_{xxxx}+(f(u))_{xx},\qquad\beta=\pm1,
\end{equation}
appeared not only in the study of the dynamics of thin inviscid
layers with free surface but also in the study of the nonlinear string, the shape-memory alloys, the propagation of waves in elastic
rods and in the continuum limit of lattice dynamics or coupled electrical circuits (see \cite{FLS} and the references therein).

Our principal aim here is to study the local well-posedness for the initial value problem associated to the sixth-order Boussinesq equation with quadratic nonlinearity \cite{cmv,daripa}:
\begin{equation}\label{sixthB}
\left\{
\begin{array}{l}
u_{tt}=u_{xx}+\beta u_{xxxx}+u_{xxxxxx}+(u^2)_{xx},\quad x\in\rr,\,\;t\geq0,\\
u(0,x)=\ff(x);\,u_t(0,x)=\psi_x(x),
\end{array}
\right.
\end{equation}
where $\beta=\pm1.$ It is worth noting that the stationary propagating localized solutions of equation \eqref{sixthB} have been investigated numerically and the
two classes of subsonic solutions corresponding to the   sign of $\beta$ have been obtained, more precisely, the monotone shapes
and the shapes with oscillatory tails \cite{cmv}.

Natural spaces to study the initial value problem above are the classical Sobolev spaces $H^s(\R)$, $s\in \R$, which are defined via the
spacial Fourier transform
\begin{equation*}
\what{g}(\xi)=\int_{\R}e^{-\ii x\xi}g(x)\;\dd x,
\end{equation*}
as the completion of the Schwarz class $\mathcal{S}(\R)$ with respect to the norm
\begin{equation*}
\|g\|_{H^{s}(\R)}=\|\langle\xi\rangle^s\;\widehat{g}(\xi)\|_{L^{2}(\R)},
\end{equation*}
where $\langle \xi \rangle = 1+|\xi|$.

Given initial data $(\phi, (\psi)_x)\in H^{s}(\R)\times H^{s-1}(\R)$ and a positive time $T>0$, we say that a function
$u: \R \times [0,T] \rightarrow \R$ is a real solution of \eqref{sixthB} if $u \in C([0,T]; H^s(\R))$ and $u$ satisfies the integral equation
\begin{equation}\label{INT}
u(t)= V_c(t)\ff+V_s(t)\psi_x+\int_{0}^{t}V_s(t-t')(u^2)_{xx}(t')\dd t'.
\end{equation}
where the two operators that constitute the free evolution are defined via Fourier transform by the formulas
\begin{gather}
V_c(t)\ff=\left( \frac{\ee^{\ii t\sqrt{{\xi}^2-\beta{\xi}^4+\xi^6}}+
\ee^{-\ii t\sqrt{{\xi}^2-\beta{\xi}^4+\xi^6}}}{2}\what{\ff}(\xi)\right)^{\vee}\\
V_s(t){\psi_x}=\left( \frac{\ee^{\ii t\sqrt{{\xi}^2-\beta{\xi}^4+\xi^6}}-
\ee^{-\ii t\sqrt{{\xi}^2-\beta{\xi}^4+\xi^6}}}{2\ii\sqrt{{\xi}^2-\beta{\xi}^4+\xi^6}}\what{\psi_x}(\xi)\right)^{\vee}.
\end{gather}

In the case that $T$ can be taken arbitrarily large, we shall say the solution is global-in-time. Here, we focus our attention only in local solutions.

Concerning the local well-posedness question, when $\beta=-1$, several results are obtained for equation \eqref{fourthB}
(so-called ``good'' Boussinesq equation) \cite{bonasachs,fanggri,farah,linares,tsut-mat}.

On the other hand  while equation \eqref{fourthB} with $\beta=1$ (so-called the ``bad'' Boussinesq equation)
only soliton type solutions are known. Moreover, taking Fourier transform, we can see that the solution of the
linearized equation $\widehat{u}$ grows as $\ee^{\pm\xi^2t}$. The same occurs for the nonlinear problem. Therefore,
to study well-posedness the component proportional to $\ee^{\xi^2t}$ has to be vanished. We refer the reader to \cite{DTT} for
results concerning this ``bad'' version using the inverse scattering approach.

The local well-posedness for dispersive equations with quadratic nonlinearities has been extensively studied in Sobolev spaces with negative
indices. The proof of these results are based on the Fourier restriction norm approach introduced by Bourgain \cite{bourgain,bourgain-2} in his study of the
nonlinear Schr\"{o}dinger (NLS) equation $\ii u_t+u_{xx}+|u|^pu=0$ and the Korteweg-de Vries (KdV) equation $u_t+u_{xxx}+uu_x=0$. This method was further developed by Kenig et al. \cite{kpv1} for the KdV Equation
and \cite{kpv2} for the quadratics nonlinear Schr\"{o}dinger equations. The original Bourgain method makes extensive use of the Strichartz
inequalities in order to derive the bilinear estimates corresponding to the nonlinearity. On the other hand, Kenig et al. simplified
Bourgain's proof and improved the bilinear estimates using only elementary techniques, such as Cauchy-Schwartz inequality and simple calculus
inequalities (see also \cite{GTV,Tao}).

In this paper, we prove local well-posedness in $H^s(\rr)$ with $s>-1/2$ for \eqref{sixthB} using the idea introduced in \cite{farah}. Indeed,
we modify the Bourgain-type space observing that the dispersion for this equation, given by the symbol
$\sqrt{{\xi}^2-\beta\xi^4+{\xi}^6}$, is in some sense related with the symbol of the KdV-type equation. This modification allow us to obtain
bilinear estimates using the same techniques as in \cite{kpv1,kpv2}.

To describe our results we define next the $X^{s,b}$ spaces related to our problem.

\begin{definition}\label{GAM}
For $s,b \in \R$, $X^{s,b}$ denotes the completion of the Schwartz class $\mathscr{S}(\R^2)$ with respect to the norm
\begin{equation}
\|u\|_{X^{s,b}}=\left\|\langle|\tau|-\gamma(\xi)\rangle^b\langle\xi\rangle^s \what{u}(\tau,\xi)\right\|_{L^{2}_{\tau,\xi}(\R^2)}
\end{equation}
where $\gamma(\xi)\equiv\sqrt{{\xi}^2-\beta{\xi}^4+\xi^6}$ and  ``$\wedge$'' denotes the time-space Fourier transform.
\end{definition}

As a consequence of this definition, we immediately have for $b> 1/2$, that $X^{s,b}$ is embedded in $C(\R;H^s(\rr))$.

We will also need the localized $X^{s,b}$ spaces defined as follows.
\begin{definition}\label{BL}
For $s,b \in \R$ and $T \geq 0$, $X^{s,b}_T$ denotes the space endowed with the norm
\begin{equation*}
\|u\|_{X^{s,b}_T}=\inf_{w\in X_{s,b}}\left\{\|w\|_{X^{s,b}}\;:\;w(t)=u(t)\,\;\mbox{on}\,\;  [0,T]\right\}.
\end{equation*}
\end{definition}
The main result of this paper reads as follows.
\begin{theorem}\label{t1.3}
Let $s>-1/2$, then for all $\ff\in H^s(\R)$ and $\psi\in H^{s-1}(\R)$, there exist
\[
T=T(\|\ff\|_{H^s(\R)},\|\psi\|_{H^{s-1}(\R)})
 \]
and a unique solution $u$ of the initial value problem associated to  equation \eqref{sixthB} with initial data $u(0)=\ff$ and $u_t(0)=\psi_x$ such that
\[
u\in C([0,T];H^s(\R))\cap X^{s,b}_T.
\]
Moreover, given $T'\in (0,T)$ there exists $R=R(T')>0$ such that giving the set
\[
W=\{(\tilde{\ff},\tilde{\psi})\in H^s(\R)\times H^{s-1}(\R):\|\tilde{\ff}-\ff\|_{H^s(\R)}^2+ \|\tilde{\psi}-\psi\|_{H^{s-1}(\R)}^2<R\}
\]
the map solution
\begin{equation*}
S:W \longrightarrow C([0,T']:H^s(\R))\cap X^{s,b}_T, \peq (\tilde{\ff},\tilde{\psi})\longmapsto u(t)
\end{equation*}
is Lipschitz. In addition, if $(\ff,\psi)\in H^{s'}(\R)\times H^{s'-1}(\R)$ with $s'>s$, then the above results hold with $s'$ instead of $s$ in the same interval $[0,T]$ with
\begin{equation*}
T=T(\|\ff\|_{H^s(\R)},\|\psi\|_{H^{s-1}(\R)}).
\end{equation*}
\end{theorem}

In some sense the previous theorem is quite surprising. There is no difference in the local theory when one considers the signs $\pm$
in front of the forth derivative term in equation \eqref{sixthB}. However, despising the sixth order term in \eqref{sixthB}, we obtain the Boussinesq equation
\eqref{fourthB}, where the ``good'' and ``bad'' models are very distinct.

We should remark that because of lack of a scaling argument for the Boussinesq-type equations, it is not clear what is the lower index $s$ where one has local well-posedness for  equation \eqref{sixthB} with initial data $u(0)=\ff$ and $u_t(0)=\psi_x$, where $(\ff,\psi)\in H^{s}(\R)\times H^{s-1}(\R)$. Here we answer, partially, this question. In fact, our main result is a negative one; it concerns in particular a kind of ill-posedness. We prove that the flow map for the Cauchy problem associated to equation \eqref{sixthB} is not smooth (more precisely $C^2$) at the origin for initial data in $H^s(\R)\times H^{s-1}(\R)$, with $s<-3$ (cf. Theorems \ref{t4.1} and \ref{t4.2}). Therefore any iterative method applied to the integral formulation of the Boussinesq equation \eqref{sixthB} always fails in this functional setting. In other words, if one can apply the contraction mapping principle to solve the integral equation corresponding to \eqref{sixthB} thus, by the implicit function Theorem, the flow-map data solution is smooth, which is a contradiction (cf. Theorem \ref{t4.2}).

Indeed our ideas are based on an argument similar to Tzvetkov \cite{tzvetkov} (see also Bourgain \cite{bourgain-2}) who established a similar result for the KdV equation. The same question was studied by Molinet, Saut and Tzvetkov \cite{MST1,MST2}, for the Benjamin-Ono (BO) equation
\begin{equation}\label{BO}
u_{t}+\mathcal{H}u_{xx}+uu_{x}=0
\end{equation}
and for the Kadomtsev-Petviashvili-I  (KPI) equation
\begin{equation}\label{KP1}
(u_{t}+uu_{x}+u_{xxx})_x-u_{yy}=0,
\end{equation}
respectively (see also \cite{aminademir}).

In all the mentioned ill-posedness results it is, in fact, proved that for a fixed $t>0$ the flow map $\ff \mapsto u(t)$ is not $C^2$ differentiable at zero. This, of course, implies that the flow map is not smooth ($C^2$) at the origin.

Unfortunately, in our case we cannot fix $t>0$ since we do not have good cancelations on the symbol $\sqrt{{\xi}^2-\beta\xi^4+{\xi}^6}$. To overcome this difficulty, we allow the variable $t$ to move. Therefore, choosing suitable characteristics functions and sending $t$ to zero to obtain our results (cf. Theorems \ref{t4.1} and \ref{t4.2}). We should remark that this kind of argument also appears in the ill-posed result of Bejenaru and Tao \cite{BT}.

The plan of this paper is as follows: in Section \ref{section2}, we state some linear estimates for the integral equation in the $X_{s,b}$
space introduced above. Bilinear estimates and the relevant counterexamples are proved in Section \ref{section3}. Finally, the ill-posedness question is
treated in Section \ref{section5}.

%%%%%%%%%%%%%%%%%%%%%%%%%%%%%%%%%%%%%%%%%%%%%%%%%%%%%%%%%%%%%%%%%%%%%%%%%%%%%%%%
%%%%%%%%%%%%%%%%%%%%%%%%%%%%%%%%%%%%%%%%%%%%%%%%%%%%%%%%%%%%%%%%%%%%%%%%%%%%%%%%
%%%%%%%%%%%%%%%%%%%%%%%%%%%%%%%%%%%%%%%%%%%%%%%%%%%%%%%%%%%%%%%%%%%%%%%%%%%%%%%%
%%%%%%%%%%%%%%%%%%%%%%%%%%%%%%%%%%%%%%%%%%%%%%%%%%%%%%%%%%%%%%%%%%%%%%%%%%%%%%%%

%%%%%%%%%%%%%%%%%%%%%%%%%%%%%%%%%%%%%%%%%%%%%%%%%%%%%%%%%%%%%%%%%%%%%%%%%%%%%%%%%%%%%%%%%%%%%%%%%%%%%%%%%%%%%%%%%%%%%%%%%%%%%%%%%%%%%%%%%%%%%%%%%%%%%%%%%%%%%%%%%%%%%%%%%%%%%%%%%%%%%%%%%%%%%%%%%%%%%%%%%%%%%%%%%%%%%%%%%%%%%%%%%%%%%%%%%%%%%%%%%%%%%%%%%%%%%%%%%%%%%%%%%%%%%%%%%%%%%%%%%%%%%%%%%%%%%%%%%%%%%%%%%%%%%%%%%%%%%%%%%%%%%%%%%%%%%%%%%%%%%%%%%%%%%%%%%%%%%%%%%%%%%%%%

\section{The Cauchy Problem}\label{section2}
Let us start this section by introducing the notation used throughout the paper. We use $c$ to denote various
constants that may vary line by line. Given any positive numbers $a$ and $b$, the notation $a \lesssim b$ means
that there exists a positive constant $c$ such that $a \leq cb$. Also, we denote $a \sim b$ when, $a \lesssim b$
and $b \lesssim a$.

Given $\theta$ be a cutoff function satisfying $\theta \in C^{\infty}_{0}(\R)$, $0\leq \theta \leq 1$, $\theta \equiv 1$ in $[-1,1]$, {\rm supp}$(\theta) \subseteq [-2,2]$ and for $0<T<1$ define $\theta_T(t)=\theta(t/T)$. In fact, to work in the $X^{s,b}$ spaces we consider another version of (\ref{INT}), that is
\begin{equation}\label{INT2}
u(t)= \theta(t)\left(V_c(t)\ff+V_s(t)\psi_x\right)+\theta_T(t)\int_{0}^{t} V_s(t-t')(u^2)_{xx}(t')\dd t'.
\end{equation}

Note that the integral equation (\ref{INT2}) is defined for all $(t,x)\in \R^2$. Moreover if $u$ is a solution of (\ref{INT2}) than $\tilde{u}=u|_{[0,T]}$ will be a solution of (\ref{INT}) in $[0,T]$.\\

In the next two lemmas, we estimate the linear and integral part of (\ref{INT2}). We refer the reader to \cite{farah} for the proofs
(see also \cite{ginibre} and \cite{GTV}).

\begin{lemma}\label{l21}
Let $u(t)$ the solution of the linear equation
\begin{eqnarray*}
\left\{
\begin{array}{l}
u_{tt}=u_{xx}+\beta u_{xxxx}+u_{xxxxxx},\\
u(0,x)=\ff(x); \peq u_t(0,x)=\psi_x(x)
\end{array} \right.
\end{eqnarray*}
with $\ff \in H^s(\rr)$ and $\psi \in H^{s-1}(\rr)$. Then there exists $c>0$ depending only on $\theta,s,b$ such that
\begin{equation}\label{LP}
\|\theta u\|_{X^{s,b}}\leq c\left(\|\ff\|_{H^s(\rr)}+\|\psi\|_{H^{s-1}(\rr)}\right).
\end{equation}
\end{lemma}
% \proof
% Taking time-space Fourier transform in $\theta(t)u(x,t)$ and setting $\gamma(\xi)=\sqrt{{\xi}^2-\beta{\xi}^4+\xi^6}$, we have
% \[\begin{split}
% (\theta(t)u(x,t))^\wedge(\tau,\xi)&=\frac{\what{\theta}(\tau-\gamma(\xi))}{2}
% \left(\what{\ff}(\xi)+\frac{\xi\what{\psi}(\xi)}{\gamma(\xi)}\right)
% +\frac{\what{\theta}(\tau+\gamma(\xi))}{2}
% \left(\what{\ff}(\xi)-\frac{\xi\what{\psi}(\xi)}{\gamma(\xi)}\right).
% \end{split}\]
% Thus, setting $h_1(\xi)=\hat{\ff}(\xi)+\frac{\xi\hat{\psi}(\xi)}{\gamma(\xi)}$ and $h_2(\xi)=\hat{\ff}(\xi)-\frac{\xi\hat{\psi}(\xi)}{\gamma(\xi)}$, we have
% \[\begin{split}
% \|\theta u\|_{X_{s,b}}^2&\leq\int_{-\infty}^{+\infty}\langle\xi\rangle^{2s}
% |h_1(\xi)|^2
% \left(\int_{-\infty}^{+\infty}\langle|\tau|-\gamma(\xi)\rangle^{2b}
% \left|\frac{\hat{\theta}(\tau-\gamma(\xi))+\hat{\theta}(\tau+\gamma(\xi))}{2}\right|^2
% \dd\tau\right)\dd\xi\\
% &\;\;+\int_{-\infty}^{+\infty}\langle\xi\rangle^{2s}
% |h_2(\xi)|^2
% \left(\int_{-\infty}^{+\infty}\langle|\tau|-\gamma(\xi)\rangle^{2b}
% \left|\frac{\hat{\theta}(\tau-\gamma(\xi))+\hat{\theta}(\tau+\gamma(\xi))}{2}\right|^2
% \dd\tau\right)\dd\xi.
% \end{split}\]
% Since $||\tau|-\gamma(\xi)|\leq \min\left\{|\tau-\gamma(\xi)|, |\tau+\gamma(\xi)|\right\}$ and $\what{\theta}$ is rapidly decreasing, we can bound the terms inside the parentheses, and the claim follows.\fim
% \indent Next we estimate the integral part of (\ref{INT2}).
\begin{lemma}\label{L22}
Let $-\frac{1}{2}<b'\leq 0\leq b \leq b'+1$ and $0<T \leq 1$ then
\begin{enumerate}[(i)]
\item   $\left\|\theta_T(t)\int_{0}^{t}g(t')\dd t'\right\|_{H^b_t} \leq T^{1-(b-b')}\|g\|_{H^{b'}_{t}}$;
\item   $\left\|\theta_T(t)\int_{0}^{t}V_s(t-t')f(u)(t')\dd t'\right\|_{X^{s,b}} \leq T^{1-(b-b')} \left\|
\left(\dfrac{\what{f(u)}(\tau,\xi)}{2\ii\gamma(\xi)}\right)^{\vee}
\right\|_{X^{s,b'}}$.
\end{enumerate}
\end{lemma}

\section{Bilinear Estimates}\label{section3}
As it is standard in the Fourier restriction method, the linear estimates given in Lemmas \ref{l21}-\ref{L22} immediately yields
Theorem \ref{t1.3} (see \cite{farah} for details) once we prove the following crucial nonlinear estimate.
\begin{theorem}\label{t1.1}
Let $s >-1/2$ and $u,v\in X^{s,-a}$. Then, there exists $c>0$ such that
\begin{equation}\label{BE}
\left\|
\left(\dfrac{|\xi|^2\what{uv}(\tau,\xi)}{2\ii\gamma(\xi)}\right)^{\vee}
\right\|_{X^{s,-a}}\leq c\left\|u\right\|_{X^{s,b}}\left\|v\right\|_{X^{s,b}},
\end{equation}
where $\vee$ denotes the inverse time-space Fourier transform, holds in the following cases
\begin{enumerate}[(i)]
\item  $s\geq 0$, $b>1/2$ and $1/6<a<1/2$,
\item $-1/2<s<0$, $b>1/2$ and $1/6<a<1/2$ such that $|s|<a$.
\end{enumerate}
Moreover, the constant $c>0$ that appears in (\ref{BE}) depends only on $a,b,s$.
\end{theorem}
%%%%%%%%%%%%%%%%%%%%%%%%%%%%%%%%%%%%%%%%
%%%%%%%%%%%%%%%%%%%%%%%%%%%%%%%%%%%%%%%%%%%%%%%%%%%%%%%%%%%%%%%%%%%%%%%%%%%%%
%%%%%%%%%%%%%%%%%%%%%%%%%%%%%%%%%%%%%%%%%%%%%%%%%%%%%%%%%%%%%%%%%%%%%%%%%%%%%
Before proceed to the proof of Theorem \ref{t1.1}, we state some elementary calculus inequalities that will be useful later.
\begin{lemma}\label{l3.1}
For $\lambda,\mu\in\rr$, $p, q>0$ and $r=\min\{ p, q, p+q-1\}$ with $p+q>1$, we have
\begin{equation}\label{CI1}
\int_\rr\dfrac{\dd x} {\langle x-\lambda\rangle^{p}\langle x-\mu\rangle^{q}}\lesssim\frac{1} {\langle\lambda-\mu\rangle^{r}}.
\end{equation}
Moreover, for $a_i\in \R$, $i=0,1,2,3$, and $q>1/3$
\begin{equation}\label{CI2}
\int_\rr\dfrac{\dd x} {\langle a_0+a_1|x|+a_2x^2+a_3|x|^3\rangle^{q}}\lesssim 1.
\end{equation}
\end{lemma}
\proof
See Lemma 4.2 in \cite{GTV} and Lemma 2.5 in \cite{BOP}.\fim
\begin{lemma}\label{l3.3}
There exists $c>0$ such that
\begin{equation}\label{LN}
\dfrac{1}{c}\leq\sup_{x,y\geq 0}\dfrac{1+\left|x-\sqrt{y^3}+\frac{\beta}{2}\sqrt{y}\right|}
{1+\left|x-\sqrt{y-\beta y^2+y^3}\right|}\leq c.
\end{equation}
\end{lemma}
\proof
Since $$\sqrt{y^3}-\frac{\beta}{2}\sqrt{y}\leq\sqrt{y-\beta y^2+y^3}\leq \sqrt{y^3}-\frac{\beta}{2}\sqrt{y}+\frac{1}{2},$$ for all $y\geq 0$ a simple computation shows the desired inequalities.
\fim
\begin{remark}\label{R1}
We should note that by using the previous lemma, we have an equivalent way to compute the $X^{s,b}$-norm, that is
\begin{equation}
\|u\|_{X^{s,b}}
\sim
\left\|\langle|\tau|-|\xi|^3+\frac{\beta}{2}|\xi|\rangle^b\langle\xi\rangle^{s} \what{u}(\tau,\xi)\right\|_{L^2_{\tau,\xi}(\rr^2)}.
\end{equation}
This equivalence will be important in the proof of Theorem \ref{t1.1}, because the symbol $\sqrt{{\xi}^2-\beta{\xi}^4+\xi^6}$ of equation \eqref{sixthB} does not have good cancelations to make use of Lemma \ref{l3.1}. Therefore, we  modify the symbols as above and work only with the algebraic relations for the KdV-type equation.
\end{remark}
%%%%%%%%%%%%%%%%%%%%%%%%%%%%%%%%%%%%%%%%%%%%%%%%%%%%%%%%%%%%%%%%%%%%%%%%%%%%%%%%
%%%%%%%%%%%%%%%%%%%%%%%%%%%%%%%%%%%%%%%%%%%%%%%%%%%%%%%%%%%%%%%%%%%%%%%%%%%%%%%%
%%%%%%%%%%%%%%%%%%%%%%%%%%%%%%%%%%%%%%%%%%%%%%%%%%%%%%%%%%%%%%%%%%%%%%%%%%%%%%%%
%%%%%%%%%%%%%%%%%%%%%%%%%%%%%%%%%%%%%%%%%%%%%%%%%%%%%%%%%%%%%%%%%%%%%%%%%%%%%%%%
%%%%%%%%%%%%%%%%%%%%%%%%%%%%%%%%%%%%%%%%%%%%%%%%%%%%%%%%%%%%%%%%%%%%%%%%%%%%%%%%
%%%%%%%%%%%%%%%%%%%%%%%%%%%%%%%%%%%%%%%%%%%%%%%%%%%%%%%%%%%%%%%%%%%%%%%%%%%%%%%%
\noindent\textbf{Proof of Theorem \ref{t1.1}.}\quad First of all observe that
\begin{equation}\label{EXI}
 \dfrac{|\xi|^2}{\sqrt{\xi^2-\beta\xi^4+\xi^6}}\leq 1, \qquad \textrm{ for all } \xi \neq 0.
\end{equation}
We will prove the theorem for the case $\beta=1$, the case $\beta=-1$ can be analogously proved.

Let $u, v\in X^{s,b}$ and define
$$
f(\tau,\xi)= \langle|\tau|-|\xi|^3+|\xi|/2\rangle^b\langle\xi\rangle^{s} \what{u}(\tau,\xi)\quad\mbox{and}\quad
g(\tau,\xi)= \langle|\tau|-|\xi|^3+|\xi|/2\rangle^b\langle\xi\rangle^{s} \what{v}(\tau,\xi).$$
Using Remark \ref{R1}, inequity \eqref{EXI} and a duality argument the desired inequality is equivalent to
\begin{equation}\label{DUA}
\left|\mathcal{W}(f,g,\ff)\right|\leq c\|f\|_{L^2_{\xi,\tau}}\|g\|_{L^2_{\xi,\tau}}\|h\|_{L^2_{\xi,\tau}},
\end{equation}
where
\begin{eqnarray*}
\mathcal{W}(f,g,h)=
\int_{\R^4}
\frac{\langle\xi\rangle^{s}
g(\tau_1,\xi_1)f(\tau_2,\xi_2) \bar{h}(\tau,\xi)}{\langle\xi_1\rangle^{s} \langle\xi_2\rangle^{s}\langle|\tau|-|\xi|^3+|\xi|/2\rangle^{a} \langle|\tau_1|-|\xi_1|^3+|\xi_1|/2\rangle^{b} \langle|\tau_2|-|\xi_2|^3+|\xi_2|/2\rangle^{b} }\dd\xi\dd\tau\dd\xi_1\dd\tau_1,
\end{eqnarray*}
where $\tau_2=\tau-\tau_1$ and $\xi_2=\xi-\xi_1$. Therefore to perform the desired estimate we need to analyze all the possible cases for the sign of $\tau$, $\tau_1$ and $\tau_2$. To do this we split $\R^4$ into the following regions:
\begin{gather*}
\Gamma_1=\left\{(\xi, \tau, \xi_1, \tau_1)\in \R^4\;:\; \tau_1, \tau_2< 0\right\},\\
\Gamma_2=\left\{(\xi, \tau, \xi_1, \tau_1)\in \R^4\;:\; \tau_1, \tau\geq 0,\; \tau_2< 0\right\},\\
\Gamma_3=\left\{(\xi, \tau, \xi_1, \tau_1)\in \R^4\;:\; \tau_1\geq 0,\; \tau_2, \tau< 0\right\},\\
\Gamma_4=\left\{(\xi, \tau, \xi_1, \tau_1)\in \R^4\;:\; \tau_1< 0,\; \tau_2, \tau\geq 0\right\},\\
\Gamma_5=\left\{(\xi, \tau, \xi_1, \tau_1)\in \R^4\;:\; \tau_1, \tau< 0,\; \tau_2\geq  0\right\},\\
\Gamma_6=\left\{(\xi, \tau, \xi_1, \tau_1)\in \R^4\;:\; \tau_1, \tau_2\geq 0\right\}.
\end{gather*}
Thus, it is suffices to prove inequality (\ref{DUA}) with $\mathcal{D}(f,g,h)$ instead of $\mathcal{W}(f,g,h)$, where
$$
\mathcal{D}(f,g,h)= \int_{\R^4} \dfrac{\langle\xi\rangle^{s}}{\langle\xi_1\rangle^{s} \langle\xi_2\rangle^{s}} \dfrac{g(\tau_1,\xi_1)f(\tau_2,\xi_2) \bar{h}(\tau,\xi)}{\langle\sigma\rangle^{a} \langle\sigma_1\rangle^{b} \langle\sigma_2\rangle^{b}}\dd\xi\dd\tau \dd\xi_1 \dd\tau_1,
$$
with $\xi_2$, $\tau_2$ and $\sigma, \sigma_1, \sigma_2$ belonging to one of the following cases
\begin{enumerate}[(I)]
\item   $\sigma=\tau+|\xi|^3-\frac{1}{2}|\xi|,\peq \sigma_1=\tau_1+|\xi_1|^3-\frac{1}{2}|\xi_1|,\peq \sigma_2=\tau_2+|\xi_2|^3-\frac{1}{2}|\xi_2|$,
\item  $\sigma=\tau-|\xi|^3+\frac{1}{2}|\xi|,\peq \sigma_1=\tau_1-|\xi_1|^3+\frac{1}{2}|\xi_1|,\peq \sigma_2=\tau_2+|\xi_2|^3-\frac{1}{2}|\xi_2|$,
\item   $\sigma=\tau+|\xi|^3-\frac{1}{2}|\xi|,\peq \sigma_1=\tau_1-|\xi_1|^3+\frac{1}{2}|\xi_1|,\peq \sigma_2=\tau_2+|\xi_2|^3-\frac{1}{2}|\xi_2|$,
\item   $\sigma=\tau-|\xi|^3+\frac{1}{2}|\xi|,\peq \sigma_1=\tau_1+|\xi_1|^3-\frac{1}{2}|\xi_1|,\peq \sigma_2=\tau_2-|\xi_2|^3+\frac{1}{2}|\xi_2|$,
\item   $\sigma=\tau+|\xi|^3-\frac{1}{2}|\xi|,\peq \sigma_1=\tau_1+|\xi_1|^3-\frac{1}{2}|\xi_1|,\peq \sigma_2=\tau_2-|\xi_2|^3+\frac{1}{2}|\xi_2|$,
\item   $\sigma=\tau-|\xi|^3+\frac{1}{2}|\xi|,\peq \sigma_1=\tau_1-|\xi_1|^3+\frac{1}{2}|\xi_1|,\peq \sigma_2=\tau_2-|\xi_2|^3+\frac{1}{2}|\xi_2|$.
\end{enumerate}
First we note that the cases $\sigma=\tau+|\xi|^3-|\xi|/2$, $\sigma_1=\tau_1-|\xi_1|^3+|\xi_2|/2$, $\sigma_2=\tau_2-|\xi_2|^3+|\xi_2|/2$ and $\sigma=\tau-|\xi|^3+|\xi|/2$, $\sigma_1=\tau_1+|\xi_1|^3-|\xi_1|/2$, $\sigma_2=\tau_2+|\xi_2|^3-|\xi_2|/2$ cannot occur, since $\tau_1< 0, \tau_2< 0$ implies $\tau<0$, and $\tau_1\geq 0$, $\tau_2\geq 0$ implies $\tau\geq 0$.
On the other hand, by applying the change of variables $(\xi, \tau, \xi_1, \tau_1)\mapsto -(\xi, \tau, \xi_1, \tau_1)$ and observing that the $L^2$-norm is preserved under the reflection operation, the cases (IV), (V), (VI) can be easily reduced, respectively, to (III), (II), (I). Moreover,  making the change of variables $\tau_2=\tau-\tau_1$, $\xi_2=\xi-\xi_1$ and then $(\xi, \tau, \xi_2, \tau_2)\mapsto -(\xi, \tau, \xi_2, \tau_2)$ the case (II) can be reduced (III). Therefore we need only establish cases (I) and (III). \\ \\

Now we first treat the inequality (\ref{DUA}) with $\mathcal{D}(f,g,h)$ in the case (I). By symmetry we can restrict ourselves to the set
\begin{equation*}
 A=\left\{(\xi, \tau, \xi_1, \tau_1)\in \R^4\;:\; |\sigma_2|\leq|\sigma_1|\right\}.
\end{equation*}
We divide $A$ into the following four subregions:
$$A_1=\{(\xi, \tau, \xi_1, \tau_1)\in A\;:\;|\xi_1|\leq10\},$$
$$A_2=\{(\xi, \tau, \xi_1, \tau_1)\in A\;:\;|\xi_1|\geq10,\;|2\xi_1-\xi|\geq|\xi_1|/2\},$$
$$A_3=\{(\xi, \tau, \xi_1, \tau_1)\in A\;:\;|\xi_1|\geq10,\;|\xi_1-\xi|\geq|\xi_1|/2,\;|\sigma_1|\leq|\sigma|\},$$
$$A_4=\{(\xi, \tau, \xi_1, \tau_1)\in A\;:\;|\xi_1|\geq10,\;|\xi_1-\xi|\geq|\xi_1|/2,\;|\sigma_1|\geq|\sigma|\}.$$
We have $A=A_1\cup A_2\cup A_3\cup A_4$. Indeed
\begin{equation*}
|\xi_1|>|2\xi_1-\xi|+|\xi_1-\xi|\geq |(2\xi_1-\xi) -(\xi_1-\xi)|=|\xi_1|.
\end{equation*}
Using the Cauchy-Schwarz and H\"older inequalities it is easy to see that
\begin{equation}\begin{split}
|Z|&\leq  \|f\|_{L^2_{\xi,\tau}(\rr^2)}\|g\|_{L^2_{\xi,\tau}(\rr^2)}\|h\|_{L^2_{\xi,\tau}(\rr^2)}
\left\|\dfrac{\langle\xi\rangle^{2s}}{\langle\sigma\rangle^{2a}} \int_{\rr^2}
\dfrac{\chi_{A_1\cup A_2\cup A_3}\;\dd\xi_1\dd\tau_1}{\langle\xi_1\rangle^{2s} \langle\xi_2\rangle^{2s} \langle\sigma_1\rangle^{2b} \langle\sigma_2\rangle^{2b}}\right\|_{L^{\infty}_{\xi,\tau}(\rr^2)}^{\frac{1}{2}}\\
&\,\,\,+\|f\|_{L^2_{\xi,\tau}(\rr^2)}\|g\|_{L^2_{\xi,\tau}(\rr^2)}\|h\|_{L^2_{\xi,\tau}(\rr^2)}
\left\|\dfrac{1}{\langle\xi_1\rangle^{2s}\langle\sigma_1\rangle^{2b}} \int_{\rr^2}\dfrac{\chi_{A_4}\langle\xi\rangle^{2s}\;\dd\xi \dd\tau}{\langle\xi_2\rangle^{2s}\langle\sigma\rangle^{2a}   \langle\sigma_2\rangle^{2b}}\right\|_{L^{\infty}_{\xi_1,\tau_1}(\rr^2)}^{\frac{1}{2}}.
\end{split}\end{equation}
Noting that $\langle\xi\rangle^{2s}\leq\langle\xi_1\rangle^{2|s|}\langle\xi_2\rangle^{2s}$, for $s\geq 0$, and $\langle\xi_2\rangle^{-2s}\leq\langle\xi_1\rangle^{2|s|}\langle\xi\rangle^{-2s}$, for $s< 0$ we have
\begin{equation}\label{Xi}
\dfrac{\langle\xi\rangle^{2s}}{\langle\xi_1\rangle^{2s} \langle\xi_2\rangle^{2s}}\leq \langle\xi_1\rangle^{\vartheta(s)}
\end{equation}
where
\begin{eqnarray*}
\vartheta(s)=
\left\{
\begin{array}{l c}
0, &\textrm{ if } s\geq 0\\
4|s|, &\textrm{ if } s\leq 0
\end{array} \right. .
\end{eqnarray*}
By employing Lemma \ref{l3.1}, it sufficient to  get bounds for
$$
J_1(\xi_1,\tau_1)=\frac{1}{\langle\sigma\rangle^{2a}}\int_{A_1\cup A_2\cup A_{3}}\frac{\langle\xi_1\rangle^{\vartheta(s)}}{\langle\tau+|\xi_2|^3-|\xi_2|/2+|\xi_1|^3-|\xi_1|/2\rangle^{2b}}\dd\xi_1
$$
and
$$
J_2(\xi,\tau)=\frac{\langle\xi_1\rangle^{\vartheta(s)}}{\langle\sigma_1\rangle^{2b}}\int_{ A_{4}}\frac{\dd\xi}{\langle\tau_1-|\xi_2|^3+|\xi_2|/2+|\xi|^3-|\xi|/2\rangle^{2a}}.
$$
\noindent\textbf{Case 1}.  Contribution of $A_1$ to $J_1$. In region $A_1$ we have $\langle\xi_1\rangle^{\vartheta(s)}\lesssim 1$. Therefore for $a>0$ and $b>1/2$ we obtain
\begin{equation*}
J_1(\xi,\tau) \lesssim \int_{|\xi_1|\leq 10}d\xi_1 \lesssim 1.
\end{equation*}
\noindent\textbf{Case 2}.  Contribution of $A_2$ to $J_1$.
In this region, we use a change of variable
$$
\eta=\tau+|\xi-\xi_1|^3-|\xi-\xi_1|/2+|\xi_1|^3-|\xi_1|/2.
$$
If $\xi_1\xi_2\geq0$, and without loss of generality $\xi_1,\xi_2\geq0$, then since $\xi\geq\xi_1\geq10$,
$$
J_1\lesssim\frac{1}{\langle\sigma\rangle^{2a}}\int\frac{\langle\xi_1\rangle^{\vartheta(s)}}
{|\xi||\xi-2\xi_1|\langle\eta\rangle^{2b}}\dd\xi_1
\lesssim
\frac{1}{\langle\sigma\rangle^{2a}}
\int\frac{\langle\xi_1\rangle^{\vartheta(s)-2}}{\langle\eta\rangle^{2b}}\dd\eta\lesssim1.
$$
If $\xi_1\xi_2\leq0$, and without loss of generality $\xi_1\geq0$ and $\xi_2\leq0$. Therefore $\xi_1\geq\xi$ and moreover
\begin{eqnarray}\label{RXI}
\xi^2+2\xi_1(\xi_1-\xi)-\frac{1}{3}=\xi_1^2+(\xi_1-\xi)^2-\frac{1}{3}\geq\frac{1}{2}\xi_1^2.
\end{eqnarray}
Using this relation, we obtain
$$
J_1\lesssim\frac{1}{\langle\sigma\rangle^{2a}}\int\frac{\langle\xi_1\rangle^{\vartheta(s)}}
{|\xi^2+2\xi_1(\xi_1-\xi)-1/3|\langle\eta\rangle^{2b}}\dd\xi_1
\lesssim
\frac{1}{\langle\sigma\rangle^{2a}}
\int\frac{\langle\xi_1\rangle^{\vartheta(s)-2}}{\langle\eta\rangle^{2b}}\dd\eta\lesssim1,
$$
for $a>0$, $b>1/2$ and $\vartheta(s)\leq2$.\\
\noindent\textbf{Case 3}.  Contribution of $A_3$ to $J_1$.
If $\xi_1\xi_2\geq0$, and without loss of generality $\xi_1,\xi_2\geq0$, then $\sigma-\sigma_1-\sigma_2=-3\xi\xi_1\xi_2$. Moreover, since $\xi\geq\xi_1\geq10$, we conclude
$$
\langle\sigma\rangle\gtrsim|\xi\xi_1\xi_2|\gtrsim|\xi_1|^3\gtrsim\langle\xi_1^3\rangle,
$$
so that Lemma \ref{l3.1} implies that
$$
J_1\lesssim\int_0^\infty\frac{\langle\xi_1\rangle^{\vartheta(s)-6a}}
{\langle\tau+(\xi-\xi_1)^3-(\xi-\xi_1)/2+\xi_1^3-\xi_1/2\rangle^{2b}}\dd\xi_1\lesssim1,
$$
for $b>1/2$ and $\vartheta(s)\leq 6a$.\\

If $\xi_1\geq0$ and $\xi_2\leq0$, then $\xi_1\geq\xi$. Therefore, by a change of variable $\eta=\tau-\xi_2^3+\xi_2/2+\xi_1^3-\xi_1/2$ and using \eqref{RXI}, we have
$$
J_1\lesssim\frac{1}{\langle\sigma\rangle^{2a}}\int\frac{\langle\xi_1\rangle^{\vartheta(s)}}
{|\xi^2+2\xi_1(\xi_1-\xi)-1/3|\langle\eta\rangle^{2b}}\dd\eta
\lesssim
\frac{1}{\langle\sigma\rangle^{2a}}
\int\frac{\langle\xi_1\rangle^{\vartheta(s)-2}}{\langle\eta\rangle^{2b}}\dd\eta\lesssim1,
$$
for $a>0$, $b>1/2$ and $\vartheta(s)\leq2$.

\noindent\textbf{Case 4}. Now we estimate $J_2(\xi_1,\tau_1)$. We use a change of variable $\eta=\tau_1+|\xi|^3-|\xi|/2-|\xi-\xi_1|^3+|\xi-\xi_1|/2$. Hence  we have $|\eta|\lesssim|\sigma_1|+|\sigma|\lesssim\langle\sigma_1\rangle$.
If $\xi,\xi_2\geq0$, then
$$|2\xi-\xi_1|=2\xi-\xi_1\geq \xi-\xi_1=|\xi-\xi_1|\geq |\xi_1|/2$$
so that
$$
J_2\lesssim
\frac{\langle\xi_1\rangle^{\vartheta(s)}}{\langle\sigma_1\rangle^{2b}}\int_{|\eta|\lesssim\langle\sigma_1\rangle}\frac{\dd\eta}
{|\xi_1||2\xi-\xi_1|\langle\eta\rangle^{2a}}
\lesssim
\frac{\langle\xi_1\rangle^{\vartheta(s)-2}}{\langle\sigma_1\rangle^{2b+2a-1}}\lesssim1,
$$
for $b>1/2$, $0<a<1/2$ and $\vartheta(s)\leq2$.\\

If $\xi\geq0$ and $\xi-\xi_1\leq0$, then one has
$$
\xi^2+(\xi-\xi_1)^2-1/3=\xi_1^2/2+(2\xi-\xi_1)^2/2 -1/3 \geq \xi_1^2/2-1/3 \geq \xi_1^2/4,
$$
so that
$$
J_2\lesssim
\frac{\langle\xi_1\rangle^{\vartheta(s)}}{\langle\sigma_1\rangle^{2b}}\int_{|\eta|\lesssim\langle\sigma_1\rangle}\frac{\dd\eta}
{|\xi^2+(\xi-\xi_1)^2-1/3|\langle\eta\rangle^{2a}}
\lesssim
\frac{\langle\xi_1\rangle^{\vartheta(s)-2}}{\langle\sigma_1\rangle^{2b+2a-1}}\lesssim1,
$$
for $b>1/2$, $0<a<1/2$ and $\vartheta(s)\leq2$.\\

Now we are going to prove the case (III). First we split $\R^4$ into the following six regions:
$$B_1=\left\{(\xi, \tau, \xi_1, \tau_1)\in \R^4\;:\;|\xi_1|\leq10\right\},$$
$$B_2=\left\{(\xi, \tau, \xi_1, \tau_1)\in \R^4\;:\;|\xi_1|\geq10,\;|\xi|\leq1\right\},$$
$$B_{3}=\left\{(\xi, \tau, \xi_1, \tau_1)\in \R^4\;:\;|\xi_1|\geq10,\;|\xi|\geq1,
;|\xi|\geq|\xi_1|/4\right\},$$
$$B_{4}=\left\{(\xi, \tau, \xi_1, \tau_1)\in \R^4\;:\;|\xi_1|\geq10,\;|\xi|\geq1,\;|\xi|\leq|\xi_1|/4,
\;\max\{|\sigma_1|,|\sigma_2|,|\sigma|\}=|\sigma|\right\},$$
$$B_{5}=\left\{(\xi, \tau, \xi_1, \tau_1)\in \R^4\;:\;|\xi_1|\geq10,\;|\xi|\geq1,\;|\xi|\leq|\xi_1|/4,
\;\max\{|\sigma_1|,|\sigma_2|,|\sigma|\}=|\sigma_1|\right\},$$
$$B_{6}=\left\{(\xi, \tau, \xi_1, \tau_1)\in \R^4\;:\;|\xi_1|\geq10,\;|\xi|\geq1,\;|\xi|\leq|\xi_1|/4,
\;\max\{|\sigma_1|,|\sigma_2|,|\sigma|\}=|\sigma_2|\right\}.$$
Using the Cauchy-Schwarz and H\"older inequalities and duality it is easy to see that
\[\begin{split}
|Z|&\leq \|f\|_{L^2_{\xi,\tau}(\rr^2)}\|g\|_{L^2_{\xi,\tau}(\rr^2)}\|h\|_{L^2_{\xi,\tau}(\rr^2)}
\left\|\dfrac{\langle\xi\rangle^{2s}}{\langle\sigma\rangle^{2a}} \int_{\rr^2}\dfrac{\chi_{B_1\cup B_3\cup B_{4}}\dd\xi_1 \dd\tau_1}{\langle\xi_1\rangle^{2s} \langle\xi_2\rangle^{2s} \langle\sigma_1\rangle^{2b} \langle\sigma_2\rangle^{2b}}\right\|_{L^{\infty}_{\xi,\tau}(\rr^2)}^{\frac{1}{2}}\\
&\,\,\,+\|f\|_{L^2_{\xi,\tau}(\rr^2)}\|g\|_{L^2_{\xi,\tau}(\rr^2)}\|h\|_{L^2_{\xi,\tau}(\rr^2)}
\left\|\dfrac{1}{\langle\xi_1\rangle^{2s}\langle\sigma_1\rangle^{2b}}
\int_{\rr^2}\dfrac{\chi_{B_2\cup B_{5}}\langle\xi\rangle^{2s}\dd\xi\dd\tau}{\langle\xi_2\rangle^{2s}\langle\sigma\rangle^{2a}   \langle\sigma_2\rangle^{2b}}\right\|_{L^{\infty}_{\xi_1,\tau_1}(\rr^2)}^{\frac{1}{2}}\\
&\,\,\,+\|f\|_{L^2_{\xi,\tau}(\rr^2)}\|g\|_{L^2_{\xi,\tau}(\rr^2)}\|h\|_{L^2_{\xi,\tau}(\rr^2)}
\left\|\dfrac{1}{\langle\xi_2\rangle^{2s}\langle\sigma_2\rangle^{2b}} \int_{\rr^2}\dfrac{\chi_{\widetilde{B}_{6}}\langle\xi_1+\xi_2\rangle^{2s}\dd\xi_1 \dd\tau_1}{\langle\xi_1\rangle^{2s}\langle\sigma_1\rangle^{2a}   \langle\sigma\rangle^{2b}}\right\|_{L^{\infty}_{\xi_2,\tau_2}(\rr^2)}^{\frac{1}{2}}.
\end{split}\]
where $\sigma$, $\sigma_1$, $\sigma_2$ were given in the condition (III) and
$$
\widetilde{B}_{6}\subset\left\{(\xi_2, \tau_2, \xi_1, \tau_1)\in \R^4\;:\;|\xi_1|\geq10,\;|\xi_1+\xi_2|\geq1\,\;|\xi_1+\xi_2|\leq|\xi_1|/4,\;
\max\{|\sigma_1|,|\sigma_2|,|\sigma|\}=|\sigma_2|\right\}.
$$
Therefore from Lemma \ref{l3.1}, it sufficient to get bounds for
$$
K_1(\xi,\tau)=\frac{1}{\langle\sigma\rangle^{2a}}\int_{B_1\cup B_3\cup B_{4}}\frac{\langle\xi_1\rangle^{\vartheta(s)}}
{\langle\tau+|\xi_2|^3 - |\xi_2|/2 -|\xi_1|^3+|\xi_1|/2\rangle^{2b}}\dd\xi_1,
$$
$$
K_2(\xi_1,\tau)=\frac{\langle\xi_1\rangle^{\vartheta(s)}}{\langle\sigma_1\rangle^{2b}}\int_{B_2\cup B_{5}}\frac{\dd\xi}
{\langle\tau_1-|\xi_2|^3 + |\xi_2|/2+|\xi|^3-|\xi|/2\rangle^{2a}},
$$
$$
K_3(\xi_2,\tau_2)=\frac{1}{\langle\sigma_2\rangle^{2b}}\int_{\widetilde{B}_{6}}\frac{\langle\xi_1\rangle^{\vartheta(s)}}
{\langle\tau_2+|\xi_1+\xi_2|^3-|\xi_1+\xi_2|/2+|\xi_1|^3-|\xi_1|/2\rangle^{2a}}\dd\xi_1,
$$
\noindent\textbf{Case 1}.  Contribution of $B_1$ to $K_1$. In region $B_1$ we have $\langle\xi_1\rangle^{\vartheta(s)}\lesssim 1$. Therefore for $a>0$ and $b>1/2$ we obtain
\begin{equation*}
K_1(\xi,\tau) \lesssim \int_{|\xi_1|\leq 10}d\xi_1 \lesssim 1.
\end{equation*}
\noindent\textbf{Case 2}.  Contribution of $B_3$ to $K_1$.
In region $B_3$, we use a change of variable $$\eta=\tau+|\xi_2|^3-|\xi_2|/2-|\xi_1|^3+|\xi_1|/2.$$ In the case $\xi_1, \xi-\xi_1\geq0$, by \eqref{RXI}, we have
\begin{equation*}
K_1\lesssim\frac{1}{\langle\sigma\rangle^{2a}}
\int\frac{\langle\xi_1\rangle^{\vartheta(s)}d\eta}{|\xi^2+2\xi_1(\xi_1-\xi)-1/3|\langle\eta\rangle^{2b}}\dd\eta
\lesssim
\frac{1}{\langle\sigma\rangle^{2a}}
\int\frac{\langle\xi_1\rangle^{\vartheta(s)-2}}{\langle\eta\rangle^{2b}}\dd\eta\lesssim1,
\end{equation*}
for $a>0$, $b>1/2$ and $\vartheta(s)\leq2$.\\

If $\xi_1\geq0$ and  $\xi-\xi_1\leq0$, then one has $|\xi-2\xi_1|=2\xi_1-\xi\geq\xi_1$, so that
$$
K_1\lesssim\frac{1}{\langle\sigma\rangle^{2a}}
\int\frac{\langle\xi_1\rangle^{\vartheta(s)}}{|\xi(\xi-2\xi_1)|\langle\eta\rangle^{2b}}\dd\eta
\lesssim
\frac{1}{\langle\sigma\rangle^{2a}}
\int\frac{\langle\xi_1\rangle^{\vartheta(s)-2}}{\langle\eta\rangle^{2b}}\dd\eta\lesssim1,
$$
for $a>0$, $b>1/2$ and $\vartheta(s)\leq2$.

\noindent\textbf{Case 3}.  Contribution of $B_4$ to $K_1$. In this region,
if $\xi_1, \xi-\xi_1\geq0$, then by definition of the set $B_{4}$, we have
$$0\leq \xi_1\leq \xi \leq \xi_1/4,$$
which is a contradiction. Therefore, this case can not happen.

Now if $\xi_1\geq0$ and $\xi-\xi_1\leq0$, then one has $|\xi-2\xi_1|=2\xi_1-\xi\geq\xi_1$.  We consider two cases. If $\xi\leq0$, then
$$
\sigma_1+\sigma_2-\sigma=3\xi\xi_1(\xi-\xi_1).
$$
If $\xi\geq0$, then
$$
\sigma_1+\sigma_2-\sigma=-2\xi^3+\xi+3\xi\xi_1(\xi-\xi_1).
$$
Since $|\xi-\xi_1|\geq 3|\xi_1|/4$, we have $|\xi|\leq|\xi\xi_1\xi_2|$ and $$2|\xi^3|\leq |\xi\xi_1(\xi-\xi_1)|/6.$$ Therefore in both cases, we have
$$
\langle\sigma\rangle\gtrsim|\xi\xi_1(\xi-\xi_1)|\gtrsim |\xi\xi_1^2|.
$$
Thus, by a change of variable $\eta=\tau-\xi^3+3\xi\xi_1(\xi-\xi_1)+\xi/2$, one gets
$$
K_1\lesssim\frac{1}{\langle\sigma\rangle^{2a}}
\int\frac{\langle\xi_1\rangle^{\vartheta(s)}}{|\xi(\xi-2\xi_1)|\langle\eta\rangle^{2b}}\dd\eta
\lesssim
\frac{1}{|\xi\xi_1^2|^{2a}}
\int\frac{\langle\xi_1\rangle^{\vartheta(s)}}{|\xi(\xi-2\xi_1)|\langle\eta\rangle^{2b}}dd\eta
\lesssim
\frac{\langle\xi_1\rangle^{\vartheta(s)-4a-1}}{|\xi|^{2a+1}}
\int\frac{\dd\eta}{\langle\eta\rangle^{2b}}
\lesssim1,
$$
for $a>1/4$, $b>1/2$ and $\vartheta(s)\leq2$.

\noindent\textbf{Case 4}.  Contribution of $B_2$ to $K_2$. First if $\xi\geq0$ and $\xi_2\leq0$, we use a change of variable $$\eta=\tau_1+\xi^3+\xi_2^3+\xi_1/2$$
and we get
$$
|\eta|\lesssim|\tau_1-\xi_1^3+\xi_1/2|+|2\xi^3-3\xi\xi_1\xi_2|\lesssim\langle\sigma_1\rangle+|\xi_1||\xi_2|
\lesssim\langle\sigma_1\rangle+\xi_1^2.
$$
Since $|\xi_1|\geq10$ and $|\xi|\leq1$ we have $|2\xi^2-\xi_1(2\xi-\xi_1)|\gtrsim |\xi_1|^ 2$. Thus we derive
\[\begin{split}
K_2&\lesssim\frac{\langle\xi_1\rangle^{\vartheta(s)}}{\langle\sigma_1\rangle^{2b}}
\int_{|\eta|\lesssim\langle\sigma_1\rangle+|\xi_1\xi_2|}\frac{\dd\eta}{|2\xi^2-\xi_1(2\xi-\xi_1)|\langle\eta\rangle^{2a}}\\
&\lesssim\frac{\langle\xi_1\rangle^{\vartheta(s)-2}}{\langle\sigma_1\rangle^{2b}}
\int_{|\eta|\lesssim\langle\sigma_1\rangle+|\xi_1|^2}\frac{\dd\eta}{\langle\eta\rangle^{2a}}
\lesssim\frac{\langle\xi_1\rangle^{\vartheta(s)-2}}{\langle\sigma_1\rangle^{2b+2a-1}}
+
\frac{\langle\xi_1\rangle^{\vartheta(s)-4a}}{\langle\sigma_1\rangle^{2b}}\lesssim1,
\end{split}\]
for $\vartheta(s)\leq\min\{2,4a\}$, $0<a<1/2$, $b>1/2$ and $2a+2b-1>0$.\\

If $\xi\geq0$ and $\xi-\xi_1\geq0$, we consider two cases. If $\xi_1\geq0$, then
$$0\leq \xi\leq \xi_1/10,$$
which is a contradiction with $\xi-\xi_1\geq0$.

So the only possible case is $\xi_1\leq0$. We use a change of variable $\eta=\tau_1-(\xi-\xi_1)^3+(\xi-\xi_1)/2+\xi^3-\xi/2$ to get
$$
|\eta|\lesssim|\tau_1+\xi_1^3|+|3\xi\xi_1(\xi-\xi_1)|\lesssim\langle\sigma_1\rangle+|\xi_1||\xi_2|
\lesssim\langle\sigma_1\rangle+\xi_1^2.
$$
Since $|\xi_1(2\xi-\xi_1)|\gtrsim |\xi_1|^ 2$, we conclude
\[\begin{split}
K_2&\lesssim\frac{\langle\xi_1\rangle^{\vartheta(s)}}{\langle\sigma_1\rangle^{2b}}
\int_{|\eta|\lesssim\langle\sigma_1\rangle+|\xi_1|^2}\frac{\dd\eta}{|\xi_1(2\xi-\xi_1)|\langle\eta\rangle^{2a}}\\
&\lesssim\frac{\langle\xi_1\rangle^{\vartheta(s)-2}}{\langle\sigma_1\rangle^{2b}}
\int_{|\eta|\lesssim\langle\sigma_1\rangle+|\xi_1|^2}\frac{\dd\eta}{\langle\eta\rangle^{2a}}
\lesssim\frac{\langle\xi_1\rangle^{\vartheta(s)-2}}{\langle\sigma_1\rangle^{2b+2a-1}}
+
\frac{\langle\xi_1\rangle^{\vartheta(s)-4a}}{\langle\sigma_1\rangle^{2b}}\lesssim1,
\end{split}\]
for $\vartheta(s)\leq\min\{2,4a\}$, $0<a<1/2$, $b>1/2$ and $2a+2b-1>0$.\\
\noindent\textbf{Case 5}.  Contribution of $B_5$ to $K_2$.
In this region, we use a change of variable $$\eta=\tau_1+|\xi|^3-|\xi|/2-|\xi_2|^3+|\xi_2|/2.$$
If $\xi,\xi_2\geq0$, then we consider two cases. If $\xi_1\geq0$, then
$$0< \xi_1\leq \xi \leq \xi_1/4,$$
which is a contradiction. Therefore, this case can not happen. If $\xi_1\leq0$, then
$$
3|\xi\xi_1\xi_2|=|\sigma_1+\sigma_2-\sigma|\lesssim\langle\sigma_1\rangle
$$
and
$$
|\eta|\lesssim|\tau_1+\xi_1^3+\xi_1/2|+|\xi\xi_1\xi_2|\lesssim\langle\sigma_1\rangle.
$$
Thus, since $|\xi_1(2\xi-\xi_1)|\gtrsim |\xi_1|^ 2$, we have
$$
K_2\lesssim\frac{\langle\xi_1\rangle^{\vartheta(s)}}{\langle\sigma_1\rangle^{2b}}
\int\frac{\dd\eta}{|\xi_1(2\xi-\xi_1)|\langle\eta\rangle^{2a}}
\lesssim\frac{\langle\xi_1\rangle^{\vartheta(s)-2}}{\langle\sigma_1\rangle^{2b}}
\int_{|\eta|\lesssim\langle\sigma_1\rangle}\frac{\dd\eta}{\langle\eta\rangle^{2a}}
\lesssim\frac{|\xi_1|^{\vartheta(s)-2}}{\langle\sigma_1\rangle^{2a+2b-1}}\lesssim1,
$$
for $\vartheta(s)\leq2$, $0<a<1/2$, $b>1/2$ and $2a+2b-1\geq0$.\\

If $\xi\geq0$ and $\xi_2\leq0$, then $\xi_1\geq0$, $\;\xi_1^2\lesssim|\xi_1(\xi_1-2\xi)|$ and
$$
\langle\sigma_1\rangle\gtrsim|\sigma_1+\sigma_2-\sigma|=|-2\xi^3+3\xi\xi_1\xi_2+\xi|\gtrsim|\xi\xi_1\xi_2|.
$$
Hence we have
$$
|\eta|\lesssim|\tau_1-\xi_1^3+\xi_1/2|+|\xi+3\xi\xi_1\xi_2|\lesssim\langle\sigma_1\rangle;
$$
and thus
$$
K_2\lesssim\frac{\langle\xi_1\rangle^{\vartheta(s)}}{\langle\sigma_1\rangle^{2b}}
\int\frac{\dd\eta}{|\xi_1(\xi_1-2\xi)-1/3|\langle\eta\rangle^{2a}}
\lesssim\frac{\langle\xi_1\rangle^{\vartheta(s)-2}}{\langle\sigma_1\rangle^{2b}}
\int_{|\eta|\lesssim\langle\sigma_1\rangle}\frac{\dd\eta}{\langle\eta\rangle^{2a}}
\lesssim\frac{|\xi_1|^{\vartheta(s)-2}}{\langle\sigma_1\rangle^{2a+2b-1}}\lesssim1,
$$
for $\vartheta(s)\leq2$, $0<a<1/2$, $b>1/2$ and $2a+2b-1\geq0$.\\
\noindent\textbf{Case 6}.  Contribution of $\widetilde{B}_6$ to $K_3$. Finally in the region $\widetilde{B}_{6}$, we have $|\xi_2|\gtrsim |\xi_1|$, therefore
$$
|\xi_1|^2\lesssim|\xi\xi_1\xi_2|\lesssim\langle\sigma_2\rangle,
$$
hence
\[\begin{split}
K_3&\lesssim\frac{1}{\langle\sigma_2\rangle^{2b}}
\int\frac{\langle\xi_1\rangle^{\vartheta(s)}}{\langle|\tau_2|+|\xi|^3-|\xi|/2+|\xi_1^3|-|\xi_1|/2\rangle^{2a}}\dd\xi_1\\
&\lesssim\langle\sigma_2\rangle^{\vartheta(s)/2-2b}
\int\frac{\dd\xi_1}{\langle|\tau_2|+|\xi|^3-|\xi|/2+|\xi_1^3|-|\xi_1|/2\rangle^{2a}}\lesssim1,
\end{split}\]
for $a>1/6$, $b>1/2$ and $\vartheta(s)\leq4b$.
\fim
%%%%%%%%%%%%%%%%%%%%%%%%%%%%%%%%%%%%%%%%%%%%%%%%%%%%%%%%%%%%%%%%%%%%%%%%%%%%%%%%
%%%%%%%%%%%%%%%%%%%%%%%%%%%%%%%%%%%%%%%%%%%%%%%%%%%%%%%%%%%%%%%%%%%%%%%%%%%%%%%%
%%%%%%%%%%%%%%%%%%%%%%%%%%%%%%%%%%%%%%%%%%%%%%%%%%%%%%%%%%%%%%%%%%%%%%%%%%%%%%%%
%%%%%%%%%%%%%%%%%%%%%%%%%%%%%%%%%%%%%%%%%%%%%%%%%%%%%%%%%%%%%%%%%%%%%%%%%%%%%%%%
%%%%%%%%%%%%%%%%%%%%%%%%%%%%%%%%%%%%%%%%%%%%%%%%%%%%%%%%%%%%%%%%%%%%%%%%%%%%%%%%
%%%%%%%%%%%%%%%%%%%%%%%%%%%%%%%%%%%%%%%%%%%%%%%%%%%%%%%%%%%%%%%%%%%%%%%%%%%%%%%%
% We finish this section with a result that will be useful in the proof of Theorem \ref{t1.3}.
% \begin{corollary}\label{c3.3}
% Let $s >-1/2$ and $a, b\in\R$ given in Theorem \ref{t1.1}. For $s'> s$ we have
% \begin{equation}\label{BE2}
% \left\|
% \left(\dfrac{|\xi|^2\what{uv}(\tau,\xi)}{2\ii\gamma(\xi)}\right)^\vee
% \right\|_{X^{s',-a}}\leq c\left\|u\right\|_{X^{s',b}}\left\|v\right\|_{X^{s,b}} + c\left\|u\right\|_{X^{s,b}}\left\|v\right\|_{X^{s',b}}.
% \end{equation}
% \end{corollary}
% \textbf{Proof. } The result is a direct consequence of Theorem \ref{t1.1} and the inequality
% \begin{equation*}
%  \langle\xi\rangle^{s'}\leq  \langle\xi\rangle^{s}\langle\xi_1\rangle^{s'-s}  +\langle\xi\rangle^{s}\langle\xi-\xi_1\rangle^{s'-s} .
% \end{equation*}
% \fim
%%%%%%%%%%%%%%%%%%%%%%%%%%%%%%%%%%%%%%%%%%%%%%%%%%%%%%%%%%%%%%%%%%%%%%%%%%%%%%%%
%%%%%%%%%%%%%%%%%%%%%%%%%%%%%%%%%%%%%%%%%%%%%%%%%%%%%%%%%%%%%%%%%%%%%%%%%%%%%%%%
%%%%%%%%%%%%%%%%%%%%%%%%%%%%%%%%%%%%%%%%%%%%%%%%%%%%%%%%%%%%%%%%%%%%%%%%%%%%%%%%
%%%%%%%%%%%%%%%%%%%%%%%%%%%%%%%%%%%%%%%%
Next we show that the bilinear estimate \eqref{BE} does not hold if $s\leq-1/2$.
More precisely,
\begin{theorem}\label{t1.4}
For any $s< -1/2$ and any $a, b\in\R$, with $a<1/2$ the estimate (\ref{BE}) fails.
\end{theorem}

The above theorem has an important consequence. It shows that our local result stated in Theorem \ref{t1.3} is sharp, in the sense that it
cannot be improved by means of the $X^{s,b}$-spaces given in Definition \ref{GAM}.\\

\noindent\textbf{Proof of Theorem \ref{t1.4}.}\quad Recall that $\gamma(\xi)=\sqrt{\xi^2-\beta\xi^4+\xi^6}$ and
let $\varrho(\xi)=\xi^3-\beta\xi/2$, $N\gg1$ and define
$$
A_N=\{(\xi,\tau)\in \R^2|N\leq \xi \leq N+N^{-\alpha}, |\tau-\varrho(\xi)|\leq 1\},
$$
where $0<\alpha<1$ will be choose later.

It is easy to see that $A_N$ contains a rectangle with $(N, 3N^2-\beta/2)$ as a vertex, with dimensions
$cN^{-2}\times N^{2-\alpha}$ and longest side pointing in the $(1,3N^2-\beta/2)$ direction.

Define $f_N(\tau,\xi)=\chi_{A_N}$ and $g_N(\tau,\xi)=\chi_{-A_N}$, then
$$
\|f_N\|_{L^2_{\tau,\xi}}\sim N^{-\alpha/2},\qquad\mbox{and}\qquad
\|g_N\|_{L^2_{\tau,\xi}}\sim N^{-\alpha/2}.
$$
Let $u_N,v_N\in X^{s,b}$ such that $$f_N(\tau,\xi)\equiv\langle|\tau|-\varrho(\xi)\rangle^b\langle\xi\rangle^s\widehat{u}_N(\tau,\xi)$$
and
$$g_N(\tau,\xi)\equiv\langle|\tau|-\varrho(\xi)\rangle^b\langle\xi\rangle^s\widehat{v}_N(\tau,\xi).$$
Therefore, from Lemma \ref{l3.3}-(\ref{LN}) and the fact that
$$
||\tau|-\varrho(\xi)|\leq\min\{|\tau-\varrho(\xi)|,|\tau+\varrho(\xi)|\},
$$
we obtain
\[\begin{split}
\left\|\left(\frac{\xi^2\widehat{u_Nv_N}(\tau,\xi)}{2\ii\gamma(\xi)}\right)^\vee\right\|_{X^{s,-a}}
&\gtrsim
\left\|\frac{\xi^2\langle\xi\rangle^s}{\gamma(\xi)\langle|\tau|-\gamma(\xi)\rangle^a}
\int_{\rr^2}\frac{f_N(\tau_1,\xi_1)g(\tau_2,\xi_2)\langle\xi_1\rangle^{-s}\langle\xi_2\rangle^{-s}\dd\tau_1\dd\xi_1}
{\langle|\tau_2|-\gamma(\xi_2)\rangle^b\langle|\tau_1|-\gamma(\xi_1)\rangle^b}\right\|_{L^2_{\tau,\xi}(\rr^2)}\\
&\gtrsim B_N,
\end{split}\]
where
$$
B_N\equiv\left\|\frac{\xi^2\langle\xi\rangle^s}{\gamma(\xi)\langle\tau-\varrho(\xi)\rangle^a}
\int_{\rr^2}\frac{f_N(\tau_1,\xi_1)g(\tau_2,\xi_2)\langle\xi_1\rangle^{-s}\langle\xi_2\rangle^{-s}\dd\tau_1\dd\xi_1}
{\langle\tau_2-\varrho(\xi_2)\rangle^b\langle\tau_1-\varrho(\xi_1)\rangle^b}\right\|_{L^2_{\tau,\xi}(\rr^2)}.
$$
From the definition of $A_N$ we have
\begin{itemize}
\item [(i)] If $(\tau_1,\xi_1)\in \textrm{supp}(f_N)$ and $(\tau_2,\xi_2)\in \textrm{supp}(g_N)$ then
\begin{equation*}
|\tau_1-\varrho(\xi_1)|\leq 1 \ppeq \textrm{ and } \ppeq |\tau_2-\varrho(\xi_2)|\leq 1.
\end{equation*}
\item [(ii)]$f\ast g(\tau,\xi)\geq \chi_{R_N}(\tau,\xi)$, where $R_N$ is the rectangle of dimensions $cN^{-2}\times{N}^{2-\alpha}$
with one of the vertices at the origin and the longest side pointing in the $(1,3N^2-\beta/2)$ direction.
\item [(iii)] $|\xi_1|\sim N$, $|\xi_2|\sim N$ and $|\xi|\leq N^{-\alpha}$.
\end{itemize}
Moreover, combining (i) and (iii) we obtain
\begin{equation}\label{TAUXI}
|\tau-\varrho(\xi)|\lesssim N^{2-\alpha}, \textrm { for all } \ppeq  |\xi|\geq N^{-\alpha}/2.
\end{equation}
Therefore (i), (ii), (iii), $(\ref{TAUXI})$ and the inequality $\xi^2/\gamma(\xi)\geq \xi$ yields
\[\begin{split}
N^{-\alpha} \peq \gtrsim\peq B_N &\gtrsim
\dfrac{N^{-(2+\alpha)s}}{N^{(2-\alpha)a}}\left\|\dfrac{|\xi|^2}{\gamma(\xi)}\chi_{R_N}\right\| _{L^2_{\tau,\xi}(\rr^2)}\\
&\gtrsim\dfrac{N^{-(2+\alpha)s}}{N^{(2-\alpha)a}}N^{-\alpha}
\left(\int_\rr\int_{\{|\xi|\geq N^{-\alpha}/2\}}{\chi}^2_{R_N}(\tau,\xi)\dd\xi\dd\tau\right)^{1/2}\\
&\gtrsim \dfrac{N^{-(2+\alpha)s}}{N^{(2-\alpha)a}}N^{-\alpha}N^{-\alpha}N^{-\alpha/2}.
\end{split}\]
Taking $N\gg1$, this inequality is possible only when
\begin{equation}\label{SA}
s\geq -\dfrac{3\alpha/2+(2-\alpha)a}{\alpha+2}.
\end{equation}

Now fix $a<1/2$ and choose $\alpha=\dfrac{1-2a}{1-a}$. Then $\alpha \in (0,1)$ and plug it into \eqref{SA} we conclude that the
estimate \eqref{BE} must fail for $s < -1/2$.
\fim

\section{Ill-posedness}\label{section5}
Before stating the main results let us define the flow-map data solution as
\begin{equation}\label{DM}
\left.
\begin{array}{c c c c}
S:&H^{s}(\R)\times H^{s-1}(\R)&\longrightarrow&C([0,T]:H^s(\R))\\
&(\ff,\psi)&\longmapsto&u(t)
\end{array} \right.
\end{equation}
where $u(t)$ is given in (\ref{INT}) below. Our ill-posedness results read as follows.
\begin{theorem}\label{t4.1}
Let $s<-3$ and any $T>0$. Then there does not exist any space $X^T$ such that

\begin{equation}\label{i}
\left\|u\right\|_{C([0,T]:H^s(\R))}\leq c\left\|u\right\|_{X^T},
\end{equation}
for all $u\in X_T$
\begin{equation}\label{ii}
\left\|V_c(t)\ff+V_s(t)\psi_x\right\|_{X^T}\leq c\left(\left\|\ff\right\|_{H^s(\R)} +\left\|\psi\right\|_{H^{s-1}(\R)}\right),
\end{equation}
for all $\ff\in H^s(\R)$, $\psi\in H^{s-1}(\R)$ and
\begin{equation}\label{iii}
\left\|\int_{0}^{t}V_s(t-t')(uv)_{xx}(t')\dd t'\right\|_{X^T}\leq c\left\|u\right\|_{X_T}\left\|v\right\|_{X^T},
\end{equation}
for all $u,v\in X^T$.
\end{theorem}
\begin{theorem}\label{t4.2}
Let $s<-3$. If there exists some $T>0$ such that the initial value problem associated to \eqref{sixthB} with initial data $u(0)=\ff$ and $u_t(0)=\psi_x$ is locally well-posed, then the flow-map data solution $S$ defined in (\ref{DM}) is not $C^2$ at zero.
\end{theorem}
\textbf{Proof of Theorem \ref{t4.1}} Suppose that there exists a space $X^T$ satisfying the conditions of the theorem for $s<-3$ and $T>0$. Let $\ff,\psi \in H^{s}(\R)$ and define $u(t)=V_c(t)\ff$, $v(t)=V_c(t)\rho$. In view of (\ref{i}), (\ref{ii}), (\ref{iii}) it is easy to see that the following inequality must hold
\begin{equation}\label{iv}
\sup_{1\leq t\leq T}\left\|\int_{0}^{t}V_s(t-t')(V_c(t')\ff V_c(t')\psi)_{xx}(t')\dd t'\right\|_{H^{s}(\R)}\leq c\left\|\ff\right\|_{H^{s}(\R)}\left\|\psi\right\|_{H^{s}(\R)}.
\end{equation}

We will see that (\ref{iv}) fails for an appropriate choice of $\ff$, $\rho$, which would lead to a contradiction.

Define
\begin{equation*}
\widehat{\ff}(\xi)=N^{-s}\chi_{[-N,-N+1]} \textrm{\para and \para} \widehat{\psi}(\xi)=N^{-s}\chi_{[N+1,N+2]},
\end{equation*}
where $\chi_A(\cdot)$ denotes the characteristic function of the set $A$.
We have
$\left\|\ff\right\|_{H^s(\R)}\sim 1$ and $\left\|\psi\right\|_{H^s(\R)} \sim 1.$ By the definitions of $V_c$, $V_s$ and Fubini's Theorem, we have
\[\begin{split}
 \left(\int_{0}^{t}V_s(t-t')(V_c(t')\ff V_c(t')\psi)_{xx}(t')\dd t'\right)^{\wedge_{x}}(\xi)&=
\int_\rr-\dfrac{|\xi|^2}{8\ii\gamma(\xi)} \widehat{\ff}(\xi_2)\widehat{\psi}(\xi_1)K(t,\xi,\xi_1)\;\dd\xi_1\\
&=\int_{A_{\xi}}-\dfrac{|\xi|^2}{8\ii\gamma(\xi)} N^{-2s}K(t,\xi,\xi_1)\;\dd\xi_1
\end{split}\]
where
\begin{equation*}
A_{\xi}=\left\{\xi_1: \xi_1\in \textrm{supp}(\widehat{\psi}) \textrm{ and } \xi_2\in \textrm{supp}(\widehat{\ff})\right\}
\end{equation*}
and
\begin{equation*}
K(t,\xi,\xi_1)\equiv \int_{0}^{t}\sin((t-t')\gamma(\xi))\cos(t'\gamma(\xi_2)) \cos(t'\gamma(\xi_1))\;\dd t'.
\end{equation*}
Note that for all $\xi_1\in \textrm{supp}(\widehat{\psi})$ and $\xi_2\in \textrm{supp}(\widehat{\ff})$ we have
\begin{equation*}
\gamma(\xi_2),\gamma(\xi_1)\sim N^3\qquad \mbox{and}\qquad 1\leq \xi \leq 3.
\end{equation*}
On the other hand, since $s<-3$, we can choose $\varepsilon>0$ such that
\begin{equation}\label{CON}
-2s-6-2\varepsilon>0.
\end{equation}
Let $t=\dfrac{1}{N^{3+\varepsilon}}$, then for $N$ sufficiently large we have
\begin{equation*}
\cos(t'\gamma(\xi_2)), \cos(t'\gamma(\xi_1))\geq 1/2
\end{equation*}
and
\begin{equation*}
\sin((t-t')\gamma(\xi))\geq c(t-t')\gamma(\xi),
\end{equation*}
for all $0\leq t'\leq t$, $1\leq \xi \leq 3$ and $\xi_1\in \textrm{supp}(\widehat{\eta})$.\\

Therefore
\begin{equation*}
K(t,\xi,\xi_1)\gtrsim \int_{0}^{t}(t-t')\gamma(\xi)\dd t'\gtrsim \gamma(\xi)\dfrac{1}{N^{6+2\varepsilon}}.
\end{equation*}

For $3/2\leq\xi\leq 5/2$ we have that $\textrm{mes}(A_{\xi})\gtrsim 1$. Thus, from (\ref{iv}) we obtain
\begin{eqnarray*}
1 &\gtrsim& \sup_{1\leq t\leq T}\left\|\int_{0}^{t}V_s(t-t')(V_c(t')\ff V_c(t')\psi)_{xx}(t')\dd t'\right\|_{H^{s}(\R)}\\
&\gtrsim& \sup_{1\leq t\leq T}\left(\int_{3/2}^{5/2}\left(1+|\xi|^2\right)^s\left| \int_{A_{\xi}}\dfrac{|\xi|^2}{8\ii\gamma(\xi)}N^{-2s} K(t,\xi,\xi_1)\dd\xi_1\right|^2\dd\xi\right)^{1/2}\\
&\gtrsim& N^{-2s-6-2\varepsilon},\quad \textrm{for all}\quad N\gg 1
\end{eqnarray*}
which is in contradiction with (\ref{CON}).
\fim
%%%%%%%%%%%%%%%%%%%%%%%%%%%%%%%%%%%%%%%%%%%%%%%%%%%%%%%%%%%%%%%%%%%%%%%%%%%%%%%%
%%%%%%%%%%%%%%%%%%%%%%%%%%%%%%%%%%%%%%%%%%%%%%%%%%%%%%%%%%%%%%%%%%%%%%%%%%%%%%%%
%%%%%%%%%%%%%%%%%%%%%%%%%%%%%%%%%%%%%%%%%%%%%%%%%%%%%%%%%%%%%%%%%%%%%%%%%%%%%%%%

\textbf{Proof of Theorem \ref{t4.2}} Let $s<-3$ and suppose that there exists $T>0$ such that the flow-map $S$ defined in (\ref{DM}) is $C^2$. When $(\ff,\psi)\in H^{s}(\R)\times H^{s-1}(\R)$, we denote by $u_{(\ff,\psi)}\equiv S(\ff,\psi)$ the solution of the IVP \eqref{sixthB} with  initial data  $u(0)=\ff$ and $u_t(0)=\psi_x$, that is
\begin{equation*}
u_{(\ff,\psi)}(t)= V_c(t)\ff+V_s(t)\psi_x+\int_{0}^{t}V_s(t-t')(u_{(\ff,\psi)}^2)_{xx}(t')\dd t'.
\end{equation*}
The Fr\'echet derivative of $S$ at $(\omega,\zeta)$ in the direction $(\ff,\bar{\ff})$ is given by
\begin{eqnarray}\label{FRECHET}
d_{(\ff,\bar{\ff})}S(\omega,\zeta)= V_c(t)\ff+V_s(t)\bar{\ff}_x+2\int_{0}^{t}V_s(t-t')(u_{(\ff,\psi)}(t') d_{(\ff,\bar{\ff})}S(\omega,\zeta)(t'))_{xx}\dd t'.
\end{eqnarray}
Using the well-posedness assumption we know that the only solution for initial data $(0,0)$ is $u_{(0,0)}\equiv S(0,0)=0$. Therefore, (\ref{FRECHET}) yields
\begin{eqnarray*}
d_{(\ff,\bar{\ff})}S(0,0)=V_c(t)\ff+V_s(t)\bar{\ff}_x.
\end{eqnarray*}
Computing the second Fr\'echet derivative at the origin in the direction $((\ff,\bar{\ff}),(\nu,\bar{\nu}))$, we obtain
\begin{eqnarray*}
d^2_{(\ff,\bar{\ff}),(\nu,\bar{\nu})}S(0,0)=2\int_{0}^{t}V_s(t-t')\left[(V_c(t')\ff+V_s(t')\bar{\ff}_x) (V_c(t')\nu+V_s(t')\bar{\nu}_x)\right]_{xx}\dd t'.
\end{eqnarray*}

Taking $\bar{\ff},\bar{\nu}=0$, the assumption of $C^2$ regularity of $S$ yields
\begin{equation*}
\sup_{1\leq t\leq T}\left\|\int_{0}^{t}V_s(t-t')(V_c(t')\ff V_c(t')\nu)_{xx}(t')\dd t'\right\|_{H^{s}(\R)}\leq c\left\|\ff\right\|_{H^{s}(\R)}\left\|\nu\right\|_{H^{s}(\R)}
\end{equation*}
which has been shown to fail in the proof of Theorem \ref{t4.1}.
\fim

%%%%%%%%%%%%%%%%%%%%%%%%%%%%%%%%%%%%%%%%%%%%%%%%%%%%%%%%%%%%%%%%%%%%%%%%%%%%%%%%
%%%%%%%%%%%%%%%%%%%%%%%%%%%%%%%%%%%%%%%%%%%%%%%%%%%%%%%%%%%%%%%%%%%%%%%%%%%%%%%%
%%%%%%%%%%%%%%%%%%%%%%%%%%%%%%%%%%%%%%%%%%%%%%%%%%%%%%%%%%%%%%%%%%%%%%%%%%%%%%%%

%\centerline{\textbf{Acknowledgment}}

%%%%%%%%%%%%%%%%%%%%%%%%%%%%%%%%%%%%%%%%%%%%%%%%%%%%%%%%%%%%%%%%%%%%%%%%%%%%%%%%
%%%%%%%%%%%%%%%%%%%%%%%%%%%%%%%%%%%%%%%%%%%%%%%%%%%%%%%%%%%%%%%%%%%%%%%%%%%%%%%%
%%%%%%%%%%%%%%%%%%%%%%%%%%%%%%%%%%%%%%%%%%%%%%%%%%%%%%%%%%%%%%%%%%%%%%%%%%%%%%%%
%%%%%%%%%%%%%%%%%%%%%%%%%%%%%%%%%%%%%%%%%%%%%%%%%%%%%%%%%%%%%%%%%%%%%%%%%%%%%%%%
%%%%%%%%%%%%%%%%%%%%%%%%%%%%%%%%%%%%%%%%%%%%%%%%%%%%%%%%%%%%%%%%%%%%%%%%%%%%%%%%

\section*{}

\end{document}